\documentclass{article}
\usepackage[centering,total={410pt,610pt}]{geometry}
\usepackage{graphicx}

\usepackage{tikz}
\usetikzlibrary{calc}
\usepackage{pgfplots}
\pgfplotsset{compat=newest}

\usepackage{enumerate}

\usepackage{amsmath}
\usepackage{amssymb}
\usepackage{amsthm}
\newtheorem{theorem}{Theorem}[section]

\newtheorem{corollary}[theorem]{Corollary}
\newtheorem{conjecture}[theorem]{Conjecture}
\theoremstyle{definition}
\newtheorem{definition}[theorem]{Definition}
\theoremstyle{remark}

\DeclareMathOperator{\prox}{prox}
\DeclareMathOperator{\argmin}{argmin}

\usepackage[ruled]{algorithm2e}

\usepackage{hyperref}
\hypersetup{colorlinks,citecolor=black,filecolor=black,linkcolor=black,urlcolor=black}

\numberwithin{equation}{section}

\usepackage{color}

\newcommand{\revise}[1]{{\color{magenta}#1}}

\title{Tight Convergence Rate in Subgradient Norm \break of the Proximal Point Algorithm}
\author{%
Guoyong Gu%
\thanks{Department of Mathematics, Nanjing University, Nanjing, 210093, China.}
\thanks{Email: ggu@nju.edu.cn. This author was supported by the NSFC grant 11671195.}
\and
Junfeng Yang\footnotemark[1]
\thanks{Email: jfyang@nju.edu.cn. This author was supported by the NSFC grant 11771208.}
}
\date{}

\begin{document}
\maketitle

\begin{abstract}
Proximal point algorithm has found many applications,
and it has been playing fundamental roles in the understanding, design, and analysis of many first-order methods.
In this paper, we derive the tight convergence rate in subgradient norm of the proximal point algorithm,
which was conjectured by Taylor, Hendrickx and Glineur [SIAM J.~Optim., 27 (2017), pp.~1283--1313].
This sort of convergence results in terms of the residual (sub)gradient norm is particularly interesting when considering dual methods,
where the dual residual gradient norm corresponds to the primal distance to feasibility.

\bigskip

\noindent\textbf{Keywords:}
proximal point algorithm,
performance estimation framework,
subgradient norm,
tight convergence rate.

\end{abstract}

\section{Introduction}
\label{Sec:Introduction}

Consider the unconstrained minimization problem
\begin{equation}
\label{model}
\min_{x\in\mathbb{R}^n} f(x).
\end{equation}
Here, the objective function $f:\mathbb{R}^n \rightarrow \mathbb{R}\cup\{+\infty\}$ is a convex, closed and proper (not necessarily differentiable) function,
which is denoted by $\mathcal{F}_{0, \infty}(\mathbb{R}^n)$ in the sequel.
Since the objective function $f$ is extended real valued, constraints such as nonnegativity, box and/or ball constraints,
can simply be absorbed into the objective function via the indicator functions.
Therefore, model (\ref{model}) encompasses a broad class of convex optimization problems.

Let $\langle \cdot, \cdot \rangle$ be the standard dot inner product 
and $B$ any symmetric positive definite matrix of order $n$. 
Assume that $\mathbb{R}^n$ is endowed with the weighted inner product $\langle x, y\rangle_B := \langle Bx, y\rangle$ for $x, y\in\mathbb{R}^n$. 
Then, the induced Euclidean norm and its dual norm are, respectively, given by
\begin{align*}
\|x\|_B = \sqrt{\langle Bx, x\rangle} \text{~~and~~}
\|u\|_{B^{-1}}&=\sqrt{\langle u, B^{-1}u\rangle}, \quad x, u\in\mathbb{R}^n.
\end{align*}
Particularly, when $B$ equals to an identity matrix, both the Euclidean norm and its dual will be denoted as $\| \cdot \|$ for simplicity.
The proximal mapping of $f\in\mathcal{F}_{0, \infty}(\mathbb{R}^n)$ is defined by
\begin{equation}\label{prox_op}
\prox_{\alpha f}(x)=\argmin_{y\in\mathbb{R}^n}\Bigl\{\alpha f(y)+\frac{1}{2}\|y-x\|_B^2 \Bigr\},\ \text{for }\alpha>0\text{~and~}x\in\mathbb{R}^n.
\end{equation}
The proximal mapping is uniquely well defined for any $x\in\mathbb{R}^n$ as the objective function in (\ref{prox_op}) is strongly convex with respect to $y$.
In this paper, we focus on the proximal point algorithm (PPA, Algotithm \ref{Alg:PPA})
for solving the aforementioned minimization problem (\ref{model}).

\begin{algorithm}[h]
\caption{Proximal point algorithm (PPA)}
\label{Alg:PPA}
\SetKwInOut{Input}{input}
\Input{Objective function: $f\in \mathcal{F}_{0,\infty}(\mathbb{R}^n)$, \break
 starting point: $x_0\in\mathbb{R}^n$,\break
 number of steps: $N$ (a positive integer),\break
 and step lengths: $\{\alpha_i\}_{i=1}^N$ with $\alpha_i>0$.}
\For {$i=1:N$}{\begin{equation}
\label{ppa}
x_i=\prox_{\alpha_i f}(x_{i-1}).
\end{equation}}
\end{algorithm}

\noindent It then follows from \eqref{prox_op}, (\ref{ppa}) and the Fermat's rule that
\begin{equation}
\label{subgradient}
0\in \alpha_i \partial f(x_i)+B(x_i-x_{i-1}), \text{ i.e., } \frac{B(x_{i-1}-x_i)}{\alpha_i}\in \partial f(x_i).
\end{equation}

PPA dates back to \cite{Mor65bsmf}.
It was firstly introduced to the optimization community in \cite{Mar70},
and later analyzed and refined by Rockafellar \cite{Roc76sicon} and G{\"u}ler \cite{Gul91sicon,Gul92siopt}.
For some recent surveys, we refer to the works of Combettes and Pesquet \cite{CP11}, Parikh and Boyd \cite{PB14fto}, and so on.

\subsection{Convergences in function and subgradient values}
The standard convergence result in terms of the function value for the PPA is provided by G{\"u}ler \cite[Theorem~2.1]{Gul91sicon}:
\begin{equation*}
f(x_N)-f(x_*)\leq \frac{R^2}{2\sum_{i=1}^{N} \alpha_i},
\end{equation*}
for any starting point $x_0\in\mathbb{R}^n$ satisfying $\|x_0-x_*\|_B\leq R$, where $R>0$ is a constant and $x_*\in\mathbb{R}^n$ is an optimal solution.
Recently, by using the performance estimation framework, this upper bound was improved by a factor of 2 by Taylor et al.:
\begin{theorem}[{\cite[Theorem~4.1]{THG17siopt}}]\label{thm:1}
Let $\{\alpha_i\}_{i\geq 1}$ be a sequence of positive step sizes
and $x_0\in\mathbb{R}^n$ some initial iterate satisfying $\|x_0-x_*\|_B\leq R$ for some optimal point $x_*$.
Any sequence $\{x_i\}_{i\geq 1}$ generated by the PPA with step sizes $\{\alpha_i\}_{i\geq 1}$
applied to a function $f\in\mathcal{F}_{0, \infty}(\mathbb{R}^n)$ satisfies
\begin{equation*}
f(x_N)-f(x_*)\leq \frac{R^2}{4\sum_{i=1}^{N} \alpha_i},
\end{equation*}
and this bound cannot be improved, even in dimension one.
\end{theorem}
\noindent The tightness of the bound provided in Theorem \ref{thm:1} is illustrated by the $l_1$-shaped one dimensional function
\begin{equation*}
f(x)=\frac{\sqrt{B}R|x|}{2\sum_{i=1}^{N} \alpha_i}
=\frac{R\|x\|_B}{2\sum_{i=1}^{N} \alpha_i}\in \mathcal{F}_{0, \infty}(\mathbb{R}),
\end{equation*}
for which $x_*=0$.
Here, $B>0$ is a constant and the starting point $x_0=-R/\sqrt{B}$ is used in the PPA.

If the residual (sub)gradient norm is used as the convergence measure,
strong numerical evidence based on performance estimation suggests the following conjecture.
\begin{conjecture}[{\cite[Conjecture~4.2]{THG17siopt}}]
\label{conj}
Let $\{\alpha_i\}_{i\geq 1}$ be a sequence of positive step sizes
and $x_0\in\mathbb{R}^n$ some initial iterate satisfying $\|x_0-x_*\|_B\leq R$ for some optimal point $x_*$.
For any sequence $\{x_i\}_{i\geq 1}$ generated by the PPA with step sizes $\{\alpha_i\}_{i\geq 1}$
on a function $f\in\mathcal{F}_{0, \infty}(\mathbb{R}^n)$,
there exists for every iterate $x_N$ a subgradient $g_N\in \partial f(x_N)$ such that
\begin{equation*}
\|g_N\|_{B^{-1}} \leq \frac{R}{\sum_{i=1}^{N} \alpha_i}.
\end{equation*}
In particular, the choice $g_N=\frac{x_{N-1}-x_N}{\alpha_N}$ is a subgradient satisfying the inequality.
\end{conjecture}
\noindent This bound cannot be improved,
as it is attained by the one-dimensional $l_1$-shaped function
\begin{equation}
\label{example}
f(x)=\frac{\sqrt{B}R |x|}{\sum_{i=1}^N\alpha_i}\in \mathcal{F}_{0, \infty}(\mathbb{R})
\end{equation}
started from $x_0=-R/\sqrt{B}$.

\subsection{The contribution and organization of the paper}
PPA has been playing fundamental roles both theoretically and algorithmically
in the optimization area \cite{Roc76sicon,Gul91sicon,Gul92siopt,YF00siopt,MS13siopt,Lie21ol}.
Even the convergence rate of the famous alternating direction method of multipliers can be established within the framework of PPA, see, e.g., ~\cite{HY12sinum}.

Convergence rate in terms of the residual (sub)gradient norm is particularly interesting when considering dual methods.
In that case, the dual residual gradient norm corresponds to the primal distance to feasibility \cite{DGN12siopt}.
Moreover, the (sub)gradient norm is more amenable than function value in nonconvex optimization.

By using performance estimation framework \cite{DT14mp, THG17siopt, THG17mp},
Conjecture~\ref{conj}, i.e., Conjecture~4.2 of \cite{THG17siopt} is proved in this paper,
which gives the first direct proof of the convergence rate in terms of the residual (sub)gradient norm for PPA.
In addition, this convergence rate for PPA is tight,
which means it is the best bound one can derive.
Note that the convergence rate in terms of the residual (sub)gradient norm can be proved
by combing the convergence rate in terms of $f(x_N)-f(x_*)$ or $\|x_N-x_*\|$ with the $L$-smoothness of $f$ \cite{THG17siopt,Nes18book},
and this kind of proof is considered to be indirect.
We want to add that the $L$-smoothness of $f$ is not required in our setting.


At first, our analysis will be restricted to the case when $B$ is the identity matrix, i.e., under the norm $\|\cdot\|$.
After that, we will show that the analysis can easily be extended to the case using a more general norm $\|\cdot\|_B$.
The rest of this paper is organized as follows.
In the next section, the performance estimation framework is briefly recalled
and the worst-case complexity bound is computed numerically via solving a semidefinite programming (SDP) reformulation.
In Section~3, the analytical optimal Lagrange multipliers is constructed based on the numerical solutions of the SDP.
In Section~4, Conjecture~\ref{conj}, i.e., Conjecture~4.2 of \cite{THG17siopt} is proved under the norm $\|\cdot\|$,
which is then extended to the case with the more general norm $\|\cdot\|_B$.
Finally, some concluding remarks are drawn in Section~5.

\section{Performance estimation and numerical results}
Performance estimation was originally developed by Drori and Teboulle \cite{DT14mp}.
Their approach is based on semidefinite relaxations,
and was taken further by Kim and Fessler \cite{KF16mp},
who derived analytically the optimized gradient method.
By using convex interpolation and smooth (strongly) convex interpolation,
the performance estimation problems can be transformed into SDP problems without any relaxation by Taylor et al.~\cite{THG17mp,THG17siopt,THG18jota}.
Recently, the performance estimation was extended by Ryu et al.~\cite{RTBG20siopt} to study operator splitting methods for monotone inclusion problems.
In \cite{GY20siopt}, the authors established tight nonergodic sublinear convergence rate of the PPA for solving maximal monotone inclusion problems.
More recently, an accelerated proximal point type algorithm for maximal monotone operator inclusion problems has been derived in \cite{Kim21mp}.

\subsection{Performance estimation and SDP reformulation}
Under the Euclidean norm $\| \cdot \|$, the worst case performance of the PPA (as given in Algorithm~\ref{Alg:PPA})
with residual subgradient norm as performance measure
can be formulated as the optimal value of the following performance estimation problem:
\begin{equation}
\label{PEP}
\tag{PEP}
\begin{aligned}
\sup_{f, x_*, x_0, x_1, \ldots, x_N, g_N} & \|g_N\|\\
\text{s.t. }\quad & f\in\mathcal{F}_{0, \infty}(\mathbb{R}^n),\\
& \|x_0-x_*\|\leq R,\\
& x_*\in\arg\min_x f(x),\\
& x_i \text{ is determined by \eqref{ppa} for } 
  1\leq i\leq N,\\
& 
g_N=\frac{x_{N-1}-x_N}{\alpha_N} \in \partial f(x_N).
\end{aligned}
\end{equation}
Suppose that $\|g_N\|$ attains some value at a $\tilde{f}\in\mathcal{F}_{0, \infty}(\mathbb{R}^n)$ with optimal solution $x_*$ and initial iterate $x_0$,
then the same value of the subgradient norm can be attained at $f(\cdot)=\tilde{f}(\cdot+x_*)-\tilde{f}(x_*)$ with optimal solution $0$ and initial iterate $x_0-x_*$.
This implies that we may assume that $x_*=0$ and $f(x_*)=0$ without affecting the optimal value of the (\ref{PEP}).


Due to the black-box property,
the PPA is determined only by the function values and the subgradients at its iterates.
According to the definition of $\mathcal{F}_{\mu,L}$-interpolation \cite[Definition~2]{THG17mp} or \cite[Definition~1.1]{THG17siopt},
the iterates of the PPA, along with the function values and subgradients at these iterates, can be considered as optimization variables instead of function $f$ itself.
As a result, the optimization problem (\ref{PEP}) can be rewritten as
\begin{align}
\notag
\sup_{f_1,\ldots, f_N, x_0, x_1,\ldots,x_N, g_1,\ldots, g_N} & \|g_N\|\\
\notag
\text{s.t. }\quad
\notag
& \|x_0\|\leq R,\\
\notag
& x_i \text{ is determined by \eqref{ppa} for } 1\leq i\leq N, \\
\notag
&  \bigl\{(0, 0, 0), (x_1, g_1, f_1), \ldots, (x_N, g_N, f_N)\bigr\} \text{ is } \mathcal{F}_{0, \infty}\text{-interpolable,}\\
\notag
&\text{where }f_i=f(x_i),\ g_i=\frac{x_{i-1}-x_i}{\alpha_i}, \text{ for } 1\leq i\leq N.
\end{align}
Remember that $(0, 0, 0)$ in the interpolation set $\bigl\{(0, 0, 0), (x_1, g_1, f_1), \ldots, (x_N, g_N, f_N)\bigr\}$ corresponds to the optimal solution $x_*=0$, $0\in \partial f(0)$, and the optimal value $f(0)=0$.
In view of \cite[Theorem~3.3]{THG17siopt},
 $\bigl\{(0, 0, 0), (x_1, g_1, f_1), \ldots, (x_N, g_N, f_N)\bigr\}$ is $\mathcal{F}_{0, \infty}$-interpolable if and only if
\begin{align*}
f_i & \geq 0+\langle 0, x_i-0\rangle,\ 1\leq i \leq N,\\
0 & \geq f_i+\langle g_i, 0-x_i\rangle,\ 1\leq i \leq N,\\
f_j & \geq f_i+\langle g_i, x_j-x_i\rangle,\ 1\leq i,j \leq N,\ i\neq j.
\end{align*}
By replacing the $\mathcal{F}_{0, \infty}$-interpolable condition with the last inequalities and eliminating $g_i$,
the optimization problem (\ref{PEP}) may further be rewritten as
%
\begin{equation}
\label{quad-PEP}
\begin{aligned}
\sup_{f_1,\ldots, f_N, x_0, x_1,\ldots,x_N} & \|x_{N-1}-x_N\|^2/\alpha_N^2\\
\text{s.t. }\quad
& \|x_0\|^2\leq R^2,\\
& f_i\geq 0,\ 1\leq i \leq N,\\
& 0\geq f_i+\bigl\langle \frac{x_{i-1}-x_i}{\alpha_i}, -x_i\bigr\rangle,\ 1\leq i \leq N,\\
& f_j\geq f_i+\bigl\langle \frac{x_{i-1}-x_i}{\alpha_i}, x_j-x_i\bigr\rangle, \ 1\leq i,j \leq N,\ i\neq j.
\end{aligned}
\end{equation}
Here, both the objective function and the constraint $\|x_0\|\leq R$ are squared.
These changes make the objective function and all the constraints simple quadratic functions over the iterates,
while neither one affects the optimal solution.

We gather the iterates of the PPA and the function values at these iterates in the following matrices
\begin{align*}
P&=\bigl[x_0, x_1, \ldots, x_N\bigr]\in \mathbb{R}^{n\times (N+1)},\\
F&=\bigl[f_1, f_2, \ldots, f_N\bigr]\in \mathbb{R}^{1\times N}.
\end{align*}
In addition, we introduce the Gram matrix  $G=P^TP\in \mathbb{S}^{N+1}_+$.
%
Notice that the rank of the Gram matrix $G$ is less than or equal to $n$.
On the other side, to reconstruct the matrix $P$ from $G$,
we need that the dimension $n$ is greater than or equal to the rank of $G$,
which is upper bounded by $N+1$.
By replacing the optimization variables with $F$ and $G$, we obtain an equivalent SDP problem, linear in $F$ and $G$, with rank constraint rank$(G) \leq n$.
Since the worst case iteration bound is considered,
the dimension $n$ can be chosen freely.
Hence, dropping the rank constraint does not incur any relaxation.


Similar reformulations have been adopted in \cite{THG17siopt,THG17mp,DT16mp,RTBG20siopt}.
Below, the SDP, with the rank constraint being dropped, is referred to as the primal SDP,
and it shares the same optimal objective function value and the Lagrange multipliers with the optimization problem (\ref{quad-PEP}).

\subsection{Numerical results}
\label{subsecNumRes}
SeDuMi \cite{Stu99oms} is utilized to solve the primal SDP, i.e., the SDP reformulation of the (\ref{PEP}).
For $R=1$, $N=20$,
and the sequence of positive step lengths $\{\alpha_i\}_{i=1}^N$ being randomly drawn from the uniform distribution between $0$ and $1$,
the absolute difference between the numerically computed optimal value of $\|g_N\|$ and the upper bound in Conjecture~\ref{conj}, namely $1/{\sum_{i=1}^{N} \alpha_i}$, is around $10^{-8}$.
In addition, the Gram matrix $G$ computed is of rank one.
The same phenomenon has been observed for many tests on other values of $N$.
This implies that there is a $f\in\mathcal{F}_{0, \infty}(\mathbb{R})$ in one-dimensional space, which attains this bound.

In fact, if the upper bound of the subgradient norm $\|g_N\|$ can be expressed as the quotient of two linear functions of $\{\alpha_i\}_{i=1}^N$,
the coefficients to these $\{\alpha_i\}_{i=1}^N$ can be found by solving a system of linear equations.
First, by multiplying both sides with the denominator,
a linear equation over the coefficients of these $\{\alpha_i\}_{i=1}^N$ can be derived
for each fixed sequence $\{\alpha_i\}_{i=1}^N$.
Then, by choosing different sequence of positive step lengths $\{\alpha_i\}_{i=1}^N$,
a system of linear equations over these coefficients can be constructed,
and from which these coefficients can be recovered numerically.



\section{The analytical Lagrange multipliers}
The optimization problem (\ref{quad-PEP}) and the primal SDP share the same Lagrange multipliers.
These Lagrange multipliers correspond to the analytical solutions of the dual SDP,
which play a key role in our proof of the conjecture.
Numerical solvers such as SeDuMi \cite{Stu99oms} only provide numerical solutions.
To find an analytical one, much more effort needs to be taken. 

\subsection{Heuristics}
To begin with, we remove the constraints that are always inactive in the primal SDP,
as the dual variables corresponding to these constraints are $0$.
For randomly chosen step lengths,
we numerically solve the primal SDP with or without certain constraints, respectively.
If the optimal values of the two  corresponding SDPs always keep the same,
then these certain constraints are considered to be inactive.
As such, our experiments indicate that the constraints corresponding to
\begin{align*}
f_i&\geq 0,\ 1\leq i \leq N-1,\\
f_j&\geq f_i+\bigl\langle \frac{x_{i-1}-x_i}{\alpha_i}, x_j-x_i\bigr\rangle, \text{ for all } |i-j|\geq 2,\ 1\leq i, j\leq N,
\end{align*}
in (\ref{quad-PEP}) can be removed from the primal SDP,
which implies that these constraints \revise{are} not supposed to play any role in the primal SDP.


Recall that the bound of the Conjecture~\ref{conj} is attained by the (one-dimensional) $l_1$-shaped function
$f(x)=\frac{R |x|}{\sum_{k=1}^N\alpha_k}$ started from $x_0=-R$.
By setting $R=1$, the iterations of the PPA (as given in Algorithm~\ref{Alg:PPA}) suggest that
\begin{align*}
P&=\Bigl[-1,\ -\frac{\alpha_{2:N}}{\alpha_{1:N}},\ -\frac{\alpha_{3:N}}{\alpha_{1:N}},\ \ldots,\ -\frac{\alpha_{N:N}}{\alpha_{1:N}},\ 0\Bigr]\in \mathbb{R}^{1\times (N+1)},\\
F&=\Bigl[\frac{\alpha_{2:N}}{(\alpha_{1:N})^2},\ \frac{\alpha_{3:N}}{(\alpha_{1:N})^2},\ \ldots,\ \frac{\alpha_{N:N}}{(\alpha_{1:N})^2},\ 0\Bigr]\in \mathbb{R}^{1\times N},
\end{align*}
from which the analytical optimal solutions for (\ref{quad-PEP}) and the primal SDP can be recovered.
Here and in the following,
to simplify the notation, we denote
\begin{equation}
\label{notationsim}
\alpha_{i:j}:=
\begin{cases}
\sum_{k=i}^j \alpha_k, & i\leq j,\\
0, & i>j.
\end{cases}
\end{equation}
The activeness of the remaining constraints in (\ref{quad-PEP}) and the primal SDP can be verified at the aforementioned optimal solution.

The Lagrange multiplier to the inequality constraint $f_N\geq 0$
in (\ref{quad-PEP}) or the corresponding constraint in the primal SDP
can be computed numerically via the method proposed in paragraph 2 of Subsection~\ref{subsecNumRes}.
Furthermore, by assuming that the multiplier is the quotient of two quadratic functions of $\{\alpha_i\}_{i=1}^N$,
the method can be generalized to compute the Lagrange multipliers corresponding to
\begin{align*}
\|x_0\|^2 &\leq 1,\ \text{(recall that } R \text{ is set to } 1 \text{)},\\
f_{N-1} &\geq f_N+\bigl\langle \frac{x_{N-1}-x_N}{\alpha_N}, x_{N-1}-x_N\bigr\rangle.
\end{align*}
However, Lagrange multipliers to the other constraints cannot be computed this way.

With the inactive constraints being removed,
there are totally $3N$ active constraints left in  the reduced  (\ref{quad-PEP}) or the primal SDP.
As a consequence, the number of Lagrangian multipliers in the Lagrangian of the reduced (\ref{quad-PEP}) is also $3N$.
On the other hand,
by setting the partial derivatives of the Lagrangian for the reduced (\ref{quad-PEP}) to $0$,
we have $2N+1$ equations over Lagrangian multipliers.
This gap implies that some Lagrangian multipliers are not unique,
which prevents us from computing them numerically.

In the Lagrangian of the reduced (\ref{quad-PEP}),
for $i=2, \ldots, N$,
by setting the partial derivative with respect to $f_i$ to $0$,
we have that the difference between the Lagrangian multipliers associated with
\begin{align*}
0 & \geq f_i+\bigl\langle \frac{x_{i-1}-x_{i}}{\alpha_{i}}, -x_{i}\bigr\rangle,\\
f_i &\geq f_{i-1}+\bigl\langle \frac{x_{i-2}-x_{i-1}}{\alpha_{i-1}}, x_i-x_{i-1}\bigr\rangle,
\end{align*}
appears in the resulting dual constraint.
Moreover, by setting the partial derivative with respect to $x_N$ to $0$,
we have the last difference,
i.e., the difference between the Lagrangian multipliers associate with
\begin{align*}
0 & \geq f_N+\bigl\langle \frac{x_{N-1}-x_{N}}{\alpha_{N}}, -x_{N}\bigr\rangle,\\
f_N &\geq f_{N-1}+\bigl\langle \frac{x_{N-2}-x_{N-1}}{\alpha_{N-1}}, x_N-x_{N-1}\bigr\rangle,
\end{align*}
equals to $(\alpha_N-\alpha_{1:N-1})/(\alpha_{1:N}\alpha_N)$.
This equation motivates us to fix one of the Lagrangian multipliers to $0$ with the other one being kept greater than or equal to $0$,
depending on the sign of $(\alpha_N-\alpha_{1:N-1})/(\alpha_{1:N}\alpha_N)$.
Similar treatments can be applied to the pair of multipliers involved in the differences for $i=2, \ldots, N-1$.
Thus, among the Lagrange multipliers of the reduced (\ref{quad-PEP}),
$N-1$ of them are fixed to $0$.

\subsection{The Lagrange multipliers}
\label{subsecmul}
The Lagrange multipliers depend on the step lengths $\{\alpha_i\}_{i=1}^N$.
In order to write down all the nonzero Lagrange multipliers,
we need to introduce a separator.
\begin{definition}[Separator]
\label{separator}
Let $\{\alpha_i\}_{i=1}^N$ be a sequence of positive step sizes,
then there exists a unique positive integer $s$
such that 
\begin{align*}
\alpha_{1:s}&>\alpha_{s+1:N},\\
\alpha_{1:s-1}&\leq \alpha_{s:N}.
\end{align*}
This positive integer $s$ is defined as the separator.
\end{definition}
\noindent Note that notation (\ref{notationsim}) is used here.
Trivially, the separator $s$ satisfies $1\leq s\leq N$ and
\begin{align*}
\alpha_{1:i}&>\alpha_{i+1:N},\ i\geq s,\\
\alpha_{1:i}&\leq \alpha_{i+1:N},\ i< s.
\end{align*}

All the nonzero Lagrange multipliers and their corresponding constraints are summarized below.
\begin{alignat*}{3}
&\frac{1}{(\alpha_{1:N})^2}\quad & &: &\quad \|x_0\|^2&\leq1,\\
&\frac{2\alpha_i}{\alpha_{i:N}\alpha_{i+1:N}} & &: &\quad 0&\geq f_i+\bigl\langle\frac{x_{i-1}-x_{i}}{\alpha_i}, -x_i\bigr\rangle,\  1\leq i \leq s-1,\\
&\frac{2\alpha_{1:i}}{\alpha_{1:N}\alpha_{i+1:N}}\quad & &: &\quad f_i&\geq f_{i+1}+\bigl\langle \frac{x_{i}-x_{i+1}}{\alpha_{i+1}}, x_i-x_{i+1}\bigr\rangle,\ 1\leq i \leq s-1,\\
&\frac{2(\alpha_{s:N}-\alpha_{1:s-1})}{\alpha_{1:N}\alpha_{s:N}}\quad & &: &\quad 0&\geq f_s+\bigl\langle \frac{x_{s-1}-x_{s}}{\alpha_s}, -x_s\bigr\rangle,\\
&\frac{2(\alpha_{1:i}-\alpha_{i+2:N})}{\alpha_{1:N}\alpha_{i+1}}\quad & &: &\quad f_i&\geq f_{i+1}+\bigl\langle \frac{x_{i}-x_{i+1}}{\alpha_{i+1}}, x_i-x_{i+1}\bigr\rangle,\ s\leq i \leq N-1,\\
&\frac{2(\alpha_{1:i}-\alpha_{i+1:N})}{\alpha_{1:N}\alpha_{i+1}}\quad & &: &\quad f_{i+1}&\geq f_{i}+\bigl\langle \frac{x_{i-1}-x_{i}}{\alpha_{i}}, x_{i+1}-x_i\bigr\rangle,\ s\leq i \leq N-1,\\
& \frac{2}{\alpha_{1:N}}\qquad & &: &\quad f_N&\geq 0.
\end{alignat*}
Here, each Lagrange multiplier and its corresponding constraint are separated by a colon.
These nonzero Lagrange multipliers together with the zero valued ones form an optimal solution to the dual SDP.
According to Definition~\ref{separator},
i.e., the definition of the separator $s$,
all the multipliers are nonnegative.

\section{Proof of the conjecture}
In a similar way as in \cite{THG17siopt,THG17mp,GY20siopt}, we prove the Conjecture~\ref{conj} in this section.
The Lagrange multipliers and their corresponding constraints summarized in the last section imply
\begin{equation}
\label{ine1}
\begin{aligned}
\Bigl\|\frac{x_{N-1}-x_{N}}{\alpha_N}\Bigr\|^2
&\leq \Bigl\|\frac{x_{N-1}-x_{N}}{\alpha_N}\Bigr\|^2-\frac{1}{(\alpha_{1:N})^2}\bigl(\|x_0\|^2-1\bigr)-A_1\\
&\quad -\frac{2(\alpha_{s:N}-\alpha_{1:s-1})}{\alpha_{1:N}\alpha_{s:N}}\Bigl(f_s+\bigl\langle \frac{x_{s-1}-x_{s}}{\alpha_s}, -x_s\bigr\rangle\Bigr)
-A_2
+\frac{2}{\alpha_{1:N}}f_N,
\end{aligned}
\end{equation}
where
\begin{align*}
A_1&=\sum_{i=1}^{s-1}\biggl[\frac{2\alpha_i}{\alpha_{i:N}\alpha_{i+1:N}}\Bigl(f_i+\bigl\langle\frac{x_{i-1}-x_{i}}{\alpha_i}, -x_i\bigr\rangle\Bigr)\\
&\qquad\qquad+\frac{2\alpha_{1:i}}{\alpha_{1:N}\alpha_{i+1:N}}\Bigl(f_{i+1}-f_{i}+\bigl\langle \frac{x_{i}-x_{i+1}}{\alpha_{i+1}}, x_i-x_{i+1}\bigr\rangle\Bigr)\biggr],\\
A_2&=\sum_{i=s}^{N-1}\biggl[\frac{2(\alpha_{1:i}-\alpha_{i+2:N})}{\alpha_{1:N}\alpha_{i+1}}\Bigl(f_{i+1}-f_i+\bigl\langle \frac{x_{i}-x_{i+1}}{\alpha_{i+1}}, x_i-x_{i+1}\bigr\rangle\Bigr)\\
&\qquad\qquad+\frac{2(\alpha_{1:i}-\alpha_{i+1:N})}{\alpha_{1:N}\alpha_{i+1}}\Bigl(f_{i}-f_{i+1}+\bigl\langle \frac{x_{i-1}-x_{i}}{\alpha_{i}}, x_{i+1}-x_i\bigr\rangle\Bigr)\biggr].
\end{align*}
By combining like terms, each $f_i$ cancels out in the right hand side of the inequality (\ref{ine1}),
and together with some rearrangement, we have
\begin{align}
\label{ine2}
\Bigl\|\frac{x_{N-1}-x_{N}}{\alpha_N}\Bigr\|^2
\leq \frac{1}{(\alpha_{1:N})^2}-A_3,
\end{align}
where
\begin{align*}
A_3
&=-\Bigl\|\frac{x_{N-1}-x_{N}}{\alpha_N}\Bigr\|^2+\frac{1}{(\alpha_{1:N})^2}\|x_0\|^2+\frac{2(\alpha_{s:N}-\alpha_{1:s-1})}{\alpha_{1:N}\alpha_{s:N}}\bigl\langle \frac{x_{s-1}-x_{s}}{\alpha_s}, -x_s\bigr\rangle\\
&\ +2\sum_{i=1}^{s-1}\Bigl[\frac{\alpha_i}{\alpha_{i:N}\alpha_{i+1:N}}\bigl\langle\frac{x_{i-1}-x_{i}}{\alpha_i}, -x_i\bigr\rangle+\frac{\alpha_{1:i}}{\alpha_{1:N}\alpha_{i+1:N}}\bigl\langle \frac{x_{i}-x_{i+1}}{\alpha_{i+1}}, x_i-x_{i+1}\bigr\rangle\Bigr]\\
&\ +2\sum_{i=s}^{N-1}\Bigl[\frac{\alpha_{1:i}-\alpha_{i+2:N}}{\alpha_{1:N}\alpha_{i+1}}\bigl\langle \frac{x_{i}-x_{i+1}}{\alpha_{i+1}}, x_i-x_{i+1}\bigr\rangle+\frac{\alpha_{1:i}-\alpha_{i+1:N}}{\alpha_{1:N}\alpha_{i+1}}\bigl\langle \frac{x_{i-1}-x_{i}}{\alpha_{i}}, x_{i+1}-x_i\bigr\rangle\Bigr].
\end{align*}

Next, we reformulate $A_3 $ as the sum of positive weighted squares.
This will imply that $A_3$ is always nonnegative.
According to the value of the separator $s$, the proof is divided into the following three categories.

\noindent (i) For $s=1$,
\begin{align*}
A_3&=\Bigl\|\frac{x_0}{\alpha_{1:N}}-\frac{\alpha_1-\alpha_{3:N}}{\alpha_1\alpha_{2}}x_1+\frac{\alpha_1-\alpha_{2:N}}{\alpha_1\alpha_{2}}x_2\Bigr\|^2
+\frac{\alpha_1-\alpha_{2:N}}{\alpha_{1:N}\alpha_1^2\alpha_2^2}\|\alpha_{3:N}x_1-\alpha_{2:N}x_2\|^2\\
&\quad +\sum_{i=s+1}^{N-1}\frac{\alpha_{1:i}-\alpha_{i+1:N}}{\alpha_{1:N}}\Bigl\|\frac{x_{i-1}}{\alpha_{i}}-\frac{\alpha_{i}+\alpha_{i+1}}{\alpha_{i}\alpha_{i+1}}x_{i}+\frac{x_{i+1}}{\alpha_{i+1}}\Bigr\|^2.
\end{align*}
\noindent (ii) For $2\leq s \leq N-1$,
\begin{align*}
A_3&=\Bigl\|\frac{x_0}{\alpha_{1:N}}-\frac{x_1}{\alpha_{2:N}}\Bigr\|^2
+\sum_{i=1}^{s-2}\frac{\alpha_{i+1}\alpha_{1:N}+2\alpha_{1:i}\alpha_{i+2:N}}{\alpha_{1:N}\alpha_{i+1}}\Bigl\|\frac{x_i}{\alpha_{i+1:N}}-\frac{x_{i+1}}{\alpha_{i+2:N}}\Bigr\|^2\\
&\quad +\frac{\alpha_{s}\alpha_{1:N}+2\alpha_{1:s-1}\alpha_{s+1:N}}{\alpha_{1:N}\alpha_s}A_4\\
&\quad +\frac{(\alpha_{1:s}-\alpha_{s+1:N})(\alpha_{s:N}-\alpha_{1:s-1})}{\alpha_{1:N}\alpha_s\alpha_{s+1}^2(\alpha_{s}\alpha_{1:N}+2\alpha_{1:s-1}\alpha_{s+1:N})}\Bigl\|\alpha_{k+2:N}x_k-\alpha_{k+1:N}x_{k+1}\Bigr\|^2\\
&\quad +\sum_{i=s+1}^{N-1}\frac{\alpha_{1:i}-\alpha_{i+1:N}}{\alpha_{1:N}}\Bigl\|\frac{x_{i-1}}{\alpha_{i}}-\frac{\alpha_{i}+\alpha_{i+1}}{\alpha_{i}\alpha_{i+1}}x_{i}+\frac{x_{i+1}}{\alpha_{i+1}}\Bigr\|^2,
\end{align*}
with
\begin{align*}
A_4=\Bigl\|\frac{x_{s-1}}{\alpha_{s:N}}-\frac{\alpha_{s+1}\alpha_{1:N}+\alpha_{s:N}(\alpha_{1:s}-\alpha_{s+1:N})}{\alpha_{s+1}(\alpha_{s}\alpha_{1:N}+2\alpha_{1:s-1}\alpha_{s+1:N})}x_s
+\frac{\alpha_{s:N}(\alpha_{1:s}-\alpha_{s+1:N})}{\alpha_{s+1}(\alpha_{s}\alpha_{1:N}+2\alpha_{1:s-1}\alpha_{s+1:N})}x_{s+1}\Bigr\|^2.
\end{align*}
\noindent (iii) For $s=N$,
\begin{align*}
A_3&=\Bigl\|\frac{x_0}{\alpha_{1:N}}-\frac{x_1}{\alpha_{2:N}}\Bigr\|^2
+\sum_{i=1}^{N-2}\frac{\alpha_{i+1}\alpha_{1:N}+2\alpha_{1:i}\alpha_{i+2:N}}{\alpha_{1:N}\alpha_{i+1}}\Bigl\|\frac{x_i}{\alpha_{i+1:N}}-\frac{x_{i+1}}{\alpha_{i+2:N}}\Bigr\|^2\\
&\quad +\frac{\alpha_N-\alpha_{1:N-1}}{\alpha_{1:N}\alpha_N^2}\|x_N\|^2.
\end{align*}
The reformulation may be verified by comparing the coefficients to $\langle x_i, x_j\rangle$ for $i,j=0, 1, \ldots, N$,
and according to Definition~\ref{separator}, the coefficients to the norm squares in the reformulations of $A_3$ are always nonnegative.

By combining the nonnegativity of $A_3$, the equation (\ref{subgradient}) and the inequality (\ref{ine2}),
we derive that for $R=1$,
\begin{equation*}
\|g_N\|^2=\Bigl\|\frac{x_{N-1}-x_{N}}{\alpha_N}\Bigr\|^2\leq \frac{1}{(\alpha_{1:N})^2}.
\end{equation*}
For $R\neq 1$, if we scale all the $\{x_i\}_{i=0, 1, \ldots, N}$ by $R$ and all the $\{f_i\}_{i=1, \ldots, N}$ by $R^2$, respectively,
we notice that the objective function value is scaled by $R^2$.
Hence, the following theorem is proved. 
\begin{theorem}
\label{theo}
Let $\{\alpha_i\}_{i\geq 1}$ be a sequence of positive step sizes
and $x_0$ some initial iterate satisfying $\|x_0-x_*\|\leq R$ for some optimal point $x_*$.
For any sequence $\{x_i\}_{i\geq 1}$ generated by the PPA with step sizes $\{\alpha_i\}_{i\geq 1}$
on any function $f\in\mathcal{F}_{0, \infty}(\mathbb{R}^n)$,
there exists for every iterate $x_N$ a subgradient $g_N\in \partial f(x_N)$ such that
\begin{equation*}
\|g_N\| \leq \frac{R}{\sum_{i=1}^{N} \alpha_i}.
\end{equation*}
In particular, the choice $g_N=\frac{x_{N-1}-x_N}{\alpha_N}$ is a subgradient satisfying the inequality.
\end{theorem}

Note that though there are some kind of heuristics in our construction of the Lagrange multipliers,
our proof of Theorem~\ref{theo} is free of these heuristics.

Up till now, our analysis is limited to the case of Euclidean norm $\| \cdot \|$.
For the more general norm $\| \cdot \|_B$,
a one to one correspondence can be formed between the feasible solutions of the Euclidean norm $\| \cdot \|$ case and the general Euclidean norm $\| \cdot \|_B$ case.
Suppose that $(f_1, \ldots, f_N, x_0, x_1, \ldots, x_N)$ is a feasible solution to (\ref{quad-PEP}),
then $\tilde{f}_1=f_1, \ldots, \tilde{f}_N=f_N$, and $\tilde{x}_0=B^{-1/2}x_0,\ \tilde{x}_1=B^{-1/2}x_1, \ldots, \tilde{x}_N=B^{-1/2}x_N$
is a feasible solution to
\begin{equation}
\label{ge-quad-PEP}
\begin{aligned}
\sup_{\tilde{f}_1, \ldots, \tilde{f}_N, \tilde{x}_0, \tilde{x}_1, \ldots, \tilde{x}_N} & \Bigl\|\frac{B(\tilde{x}_{N-1}-\tilde{x}_N)}{\alpha_N}\Bigr\|_{B^{-1}}^2 \\
\text{s.t. }\quad
& \|\tilde{x}_0\|_B^2\leq R^2,\\
& \tilde{f}_i\geq 0,\ 1\leq i \leq N,\\
& 0\geq \tilde{f}_i+\Bigl\langle \frac{B(\tilde{x}_{i-1}-\tilde{x}_i)}{\alpha_i}, -\tilde{x}_i\Bigr\rangle,\ 1\leq i \leq N,\\
& \tilde{f}_j\geq \tilde{f}_i+\Bigl\langle \frac{B(\tilde{x}_{i-1}-\tilde{x}_i)}{\alpha_i}, \tilde{x}_j-\tilde{x}_i\Bigr\rangle, \ 1\leq i,j \leq N,\ i\neq j,
\end{aligned}
\end{equation}
and vice versa.
The optimization problem (\ref{ge-quad-PEP}) is just a reformulation of the performance estimation problem under the general norm $\| \cdot \|_B$.
Notice that the two optimization problems (\ref{quad-PEP}) and (\ref{ge-quad-PEP}) have the same objective function value
at these one to one correspondent feasible solutions.
As a consequence, they also have the same optimal value.
Hence, we have the following corollary.
\begin{corollary}[{\cite[Conjecture~4.2]{THG17siopt}}]
Let $\{\alpha_i\}_{i\geq 1}$ be a sequence of positive step sizes
and $x_0$ some initial iterate satisfying $\|x_0-x_*\|_B\leq R$ for some optimal point $x_*$.
For any sequence $\{x_i\}_{i\geq 1}$ generated by the PPA (based on the general norm $\|\cdot\|_B$) with step sizes $\{\alpha_i\}_{i\geq 1}$
on a function $f\in\mathcal{F}_{0, \infty}(\mathbb{R}^n)$,
there exists for every iterate $x_N$ a subgradient $g_N\in \partial f(x_N)$ such that
\begin{equation*}
\|g_N\|_{B^{-1}} \leq \frac{R}{\sum_{i=1}^{N} \alpha_i}.
\end{equation*}
In particular, the choice $g_N=\frac{B(x_{N-1}-x_N)}{\alpha_N}$ is a subgradient satisfying the inequality.
\end{corollary}
\noindent Note that this bound cannot be improved,
as it is attained on the (one dimensional) $l_1$-shaped function (\ref{example}).

\section{Concluding remarks}
PPA has been playing fundamental roles both theoretically and algorithmically in the optimization area.
The convergence results in terms of the residual (sub)gradient norm are particularly interesting when considering dual methods \cite{DGN12siopt}.
Moreover, the (sub)gradient norm is more amenable than function value in nonconvex optimization.

By using performance estimation framework, Conjecture~4.2 in \cite{THG17siopt} is proved,
which gives the first direct proof of the convergence rate in terms of the residual (sub)gradient norm for PPA.
To the best of our knowledge, even for PPA with universal step length, i.e., freezing all step lengths at the same positive constant,
a direct proof of the convergence rate in terms of the residual (sub)gradient norm is not known before.
In addition, this convergence rate is tight, which means it is the best bound one can derive.


\end{document}